# SUPPLEMENT TO "ERRATUM: HIGHER ORDER ELICITABILITY AND OSBAND'S PRINCIPLE"

By Tobias Fissler[*] and Johanna F. Ziegel

*WU Vienna University of Economics and Business and University of Bern*

**Abstract.** This note corrects conditions in Proposition 3.4 and Theorem 5.2(ii) and comments on imprecisions in Propositions 4.2 and 4.4 in Fissler and Ziegel (2016a).

**1. Proposition 3.4.** As detailed in Brehmer (2017, Subsections 1.3.2 and 1.3.3) there are two technicalities that need to be resolved in Fissler and Ziegel (2016a, Proposition 3.4): Firstly, due to the particular choice of the integration path in the original version of Fissler and Ziegel (2016a, Proposition 3.4), the image of the integration path is not necessarily contained in int(A). Secondly, one needs to assume that the identification function $V$ is locally bounded *jointly* in the two components. Proposition 1 gives a corrected version of Fissler and Ziegel (2016a, Proposition 3.4).

PROPOSITION 1. *Assume that* $\mathrm{int}(\mathsf{A}) \subseteq \mathbb{R}^k$ *is simply connected and let* $T\colon \mathcal{F} \to \mathsf{A}$ *be a surjective, elicitable and identifiable functional with a strict $\mathcal{F}$-identification function* $V\colon \mathsf{A} \times \mathsf{O} \to \mathbb{R}^k$ *and a strictly $\mathcal{F}$-consistent scoring function* $S\colon \mathsf{A} \times \mathsf{O} \to \mathbb{R}$. *Suppose that Assumption (V1), (V2), (S1) from Fissler and Ziegel (2016a) are satisfied. Let $h$ be the matrix-valued function appearing at Fissler and Ziegel (2016a, Equation (3.2)).*

*For any $F \in \mathcal{F}$ and any points $x, z \in \mathrm{int}(\mathsf{A})$ and any integration path $\gamma\colon [0,1] \to \mathrm{int}(\mathsf{A})$ with $\gamma(0) = z$, $\gamma(1) = x$ the score difference is necessarily of the form*

$$(1) \quad \bar{S}(x, F) - \bar{S}(z, F) = \int_\gamma \mathrm{d}\bar{S}(\cdot, F) = \int_0^1 \left( h(\gamma(\lambda)) \bar{V}(\gamma(\lambda), F) \right)^\top \gamma'(\lambda)\, \mathrm{d}\lambda.$$

*Moreover, if Assumptions (F1) and (VS1) from Fissler and Ziegel (2016a) are satisfied and $V$ is locally bounded, then there is a Lebesgue null set $N \subseteq$*

---

[*]Tobias Fissler gratefully received funding of his Chapman Fellowship by the Department of Mathematics at Imperial College London.

*MSC 2010 subject classifications:* Primary 62C99; secondary 91B06

$\mathsf{A} \times \mathsf{O}$ *such that for all* $(x,y) \in N^c$ *and for all* $(z,y) \in N^c$ *it necessarily holds that*

$$(2) \quad S(x,y) - S(z,y) = \int_\gamma \mathrm{d}S(\cdot,y) = \int_0^1 \Big(h\big(\gamma(\lambda)\big)V\big(\gamma(\lambda),y\big)\Big)^\top \gamma'(\lambda)\,\mathrm{d}\lambda\,,$$

*where* $\gamma\colon [0,1] \to \mathrm{int}(\mathsf{A})$ *is an integration path with* $\gamma(0) = x$, $\gamma(1) = z$.

PROOF. Equation (1) follows from Fissler and Ziegel (2016a, Theorem 3.2) and Königsberger (2004, Satz 2, p. 183). The proof of (2) follows the lines of the original proof in Fissler and Ziegel (2016b); cf. Brehmer (2017, Theorem 1.31) for details. □

We would like to briefly discuss the implications of Proposition 1 on other results in Fissler and Ziegel (2016a) where the original Proposition 3.4 is used.

Proposition 4.2(ii): Concerning the conditions, one needs to add that $V$ is locally bounded. On the other hand, it is sufficient to impose that $\mathrm{int}(\mathsf{A})$ is simply connected and we can dispense with the assumption of a star domain.

The proof of the proposition on p. 1701 of the original paper basically remains the same. The sentence around the second displayed equation should read as follows: For necessity, we apply Proposition 1 and part (i) such that for almost all $(x,y), (z,y) \in \mathsf{A} \times \mathsf{O}$ there is an integration path $\gamma\colon [0,1] \to \mathrm{int}(\mathsf{A})$ with $\gamma(0) = z$ and $\gamma(1) = x$ such that

$$\begin{aligned}
S(x,y) - S(z,y) &= \int_0^1 \Big(h\big(\gamma(\lambda)\big)V\big(\gamma(\lambda),y\big)\Big)^\top \gamma'(\lambda)\,\mathrm{d}\lambda \\
&= \sum_{m=1}^k \int_0^1 g_m\big(\gamma(\lambda)\big)V_m\big(\gamma(\lambda),y\big)\gamma'_m(\lambda)\,\mathrm{d}\lambda \\
&= \sum_{m=1}^k \int_{z_k}^{x_k} g_m(v)V_m(v,y)\,\mathrm{d}v.
\end{aligned}$$

The rest remains unchanged.

Proposition 4.4(ii): Again, we need to impose that $V$ is locally bounded. On the other hand, assuming that $\mathrm{int}(\mathsf{A})$ is simply connected is sufficient. Note that the specific identification function $V$ is an affine function in $x$ such that we do not need to explicitly assume (V3). What we need, however, on top of Assumption (S2) is that the expected score $\bar{S}(\cdot, F)$ is twice continuously differentiable on $\mathrm{int}(\mathsf{A})$ for all $F \in \mathcal{F}$.

This additional smoothness implies that the function $h$ is continuously differentiable and that the first identity in (4.4) holds on $\text{int}(\mathsf{A})$. Then, a double application of Königsberger (2004, Satz 6, p. 193) yields that there exists a three times differentiable function $\phi\colon \text{int}(\mathsf{A}) \to \mathbb{R}$ with first order partial derivatives $\partial_m \phi$ and Hessian $h$. For $z \in \text{int}(\mathsf{A})$ and $x \in \text{int}(\mathsf{A})$, let $\gamma\colon [0,1] \to \text{int}(\mathsf{A})$ be an integration path such that $\gamma(0) = z$ and $\gamma(1) = x$. Then, partial integration yields

$$S(x,y) - S(z,y) = \int_\gamma \mathrm{d}S(\cdot, y) = \int_0^1 \Big(h\big(\gamma(\lambda)\big)V\big(\gamma(\lambda), y\big)\Big)^\top \gamma'(\lambda)\, \mathrm{d}\lambda$$

$$= \sum_{m=1}^k \int_0^1 \sum_{i=1}^k \gamma_i'(\lambda) h_{im}\big(\gamma(\lambda)\big) V_m\big(\gamma(\lambda), y\big)\, \mathrm{d}\lambda$$

$$(3) \quad = \sum_{m=1}^k \partial_m \phi\big(\gamma(\lambda)\big) V_m\big(\gamma(\lambda), y\big)\Big|_0^1 - \int_0^1 \partial_m \phi\big(\gamma(\lambda)\big) q(y) \gamma_m'(\lambda)\, \mathrm{d}\lambda$$

$$(4) \quad = \sum_{m=1}^k \partial_m \phi(x) \big(q(y) x_m - p_m(y)\big) - \phi(x) q(y)$$
$$\quad - \partial_m \phi(z)\big(q(y) z_m - p_m(y)\big) + \phi(z) q(y).$$

Notice that this does not depend on the particular choice of $\phi$. Indeed, any other function $\tilde{\phi}$ with the same properties can be written as $\tilde{\phi}(x) = \phi(x) + \beta^\top x + \alpha$ for some $\alpha \in \mathbb{R}$, $\beta \in \mathbb{R}^k$. Then one can easily calculate that using $\tilde{\phi}$ instead of $\phi$ in (3) does not affect the result in (4).

Theorem 5.2(iii) and implicitly Corollary 5.5(iii): Here, the identification function given at Equation (5.4) is locally bounded, such that we do not need to impose this condition. Moreover, it is again enough to assume that $\text{int}(\mathsf{A})$ is simply connected. Note that this implies that the projection of $\text{int}(\mathsf{A})$ onto the $k$th component is therefore an interval. Therefore, we do not need any involved argument for the existence of $\mathcal{G}_k$. The rest relies again on a partial integration argument similar to the ones above.

**2. Theorem 5.2(ii).** The error in the proof of Theorem 5.2(ii) in Fissler and Ziegel (2016a) can be found on p. 1702 directly under the equation defining the term $R_2$. We wrote "Due to the assumptions, the term $G_r(y) + (p_r/q_r)G_k(w)y$ is increasing in $y \in [t_r, x_r]$" where $w := \min(x_k, t_k)$. However, our assumptions do not ensure that the interval $[t_r, x_r]$ is necessarily contained in the set

$$\mathsf{A}'_{r,w} := \{x_r \colon \exists (z_1, \ldots, z_k) \in \mathsf{A},\ x_r = z_r,\ w = z_k\}.$$

Note that $\mathsf{A}'_{r,w}$ is the projection of the $w$-section of $\mathsf{A}$ to the $r$th component while $\mathsf{A}'_r := \{x_r \colon \exists (z_1, \ldots, z_k) \in \mathsf{A},\ x_r = z_r\}$ is simply the projection of $\mathsf{A}$ to the $r$th component. Hence, the condition at Equation (5.3) cannot readily be applied.

In Subsection 2.1 we give a counterexample which demonstrates that the complication described above can indeed lead to a scoring function satisfying the conditions of Theorem 5.2(ii) in Fissler and Ziegel (2016a) which is not consistent. In Subsection 2.2 we introduce a condition on the action domain which, in combination with the conditions in Theorem 5.2(ii), ensures (strict) consistency of the scoring function. We end with some remarks as to when this additional condition is satisfied.

2.1. *Counterexample.* Using the same notation as in Fissler and Ziegel (2016a), we confine ourselves to presenting a counterexample for $k = 2$ and $\alpha \in (0, 1/2)$ (counterexamples for $\alpha \in [1/2, 1)$ can be constructed in a similar manner). Consider the convex action domain $\mathsf{A} = \{(x_1, x_2) \in \mathbb{R}^2 \colon x_1 \geq 0,\ |x_2| \leq x_1\}$. Let $\mathcal{G}_2 \colon \mathbb{R} \to \mathbb{R}$ be a strictly convex twice differentiable function and define $G_1 \colon [0, \infty) \to \mathbb{R}$ via

$$G_1(s) = \mathcal{G}_2(-s)/\alpha, \quad s \geq 0.$$

Observe that the condition at (5.3) in Fissler and Ziegel (2016a) holds since for $(x_1, x_2) \in \mathsf{A}$

$$G'_1(x_1) + \mathcal{G}'_2(x_2)/\alpha = \left(\mathcal{G}'_2(x_2) - \mathcal{G}'_2(-x_1)\right)/\alpha \begin{cases} = 0, & \text{if } x_2 = -x_1, \\ > 0, & \text{if } x_2 > -x_1, \end{cases}$$

due to the fact that $\mathcal{G}_2$ is strictly convex. Consider the specific choice $\mathcal{G}_2 = \exp$ and $a \equiv 0$. This results in the score

$$\begin{aligned} S(x_1, x_2, y) &= (\mathbb{1}\{y \leq x_1\} - \alpha) \exp(-x_1)/\alpha - \mathbb{1}\{y \leq x_1\} \exp(-y)/\alpha \\ &\quad + \exp(x_2)\big(x_2 + (\mathbb{1}\{y \leq x_1\} - \alpha)x_1/\alpha - \mathbb{1}\{y \leq x_1\}y/\alpha\big) - \exp(x_2). \end{aligned}$$

To demonstrate that the score fails to be consistent for the pair $T = (\mathrm{VaR}_\alpha, \mathrm{ES}_\alpha)$ with $\alpha = 0.05$ consider the following distributions: First, $F$ a point-distribution in 0. Then $T(F) = (0, 0)$ and we obtain $\bar{S}(T(F), F) = S(0, 0, 0) = -2$ and, for example, $S(2, -1.8, 0) \approx -11.61$. Second, $F$ a normal distribution with mean $\mu = 0.2$ and standard deviation $\sigma = 0.1$. Then $T(F) \approx (0.0355, -0.0063)$. A numerical integration yields that $\bar{S}(T(F), F) \approx -5.36$ and, for example, $\bar{S}(2, -1.8, F) \approx -8.76$.

2.2. *A sufficient condition for Theorem 5.2(ii).*

PROPOSITION 2. *Let the conditions of Theorem 5.2 in Fissler and Ziegel (2016a)(ii) hold. Moreover, assume that for any $x \in \mathsf{A}$ and any $t \in T(\mathcal{F})$, there exists a finite sequence $(z^{(n)})_{n=0,\ldots,N} \subseteq \mathsf{A}$ such that*

1. $z^{(0)} = x$;
2. $t_r, z_r^{(N)} \in \mathsf{A}'_{r,\min\{t_k, z_k^{(N)}\}}$ *for all $r \in \{1, \ldots, k-1\}$;*
3. *for all $n \in \{1, \ldots, N\}$ there is an index $r \in \{1, \ldots, k\}$ such that*
    a) $t_r \leq z_r^{(n)} \leq z_r^{(n-1)}$ *or* $z_r^{(n-1)} \leq z_r^{(n)} \leq t_r$; *and*
    b) $z_m^{(n)} = z_m^{(n-1)}$ *for all $m \in \{1, \ldots, k\} \setminus \{r\}$;*
4. *for all $n \in \{1, \ldots, N\}$ it holds that if $x_k < t_k$ and $z_k^{(n)} \neq z_k^{(n-1)}$, then*
$$z_k^{(n-1)} < z_k^{(n)} \leq -B(z^{(n)}, t) = -B(z^{(n-1)}, t),$$

*where*

$$(5) \quad B(x,t) := -t_k + \sum_{m=1}^{k-1} \frac{p_m}{q_m}(t_m - x_m)(q_m - \mathbb{1}\{t_m < x_m\}) \quad x, t \in \mathbb{R}^k.$$

*Then the scoring function $S$ defined at (5.2) is $\mathcal{F}$-consistent for $T$. If additionally, the distributions in $\mathcal{F}$ have unique $q_m$-quantiles, $m \in \{1, \ldots, k-1\}$, $\mathcal{G}_k$ is strictly convex and the functions given at (5.3) are strictly increasing, then $S$ is strictly $\mathcal{F}$-consistent for $T$.*

First, we provide some intuition how the proof works, and comment on the additional condition in Proposition 2 concerning the existence of the sequence $(z^{(n)})_{n=0,\ldots,N} \subseteq \mathsf{A}$. The idea of the proof is to show that for any $F \in \mathcal{F}$, $t = T(F)$, $x \in \mathsf{A}$, there is a sequence $(z^{(n)})_{n=0,\ldots,N} \subseteq \mathsf{A}$, which basically constitutes a path in $\mathsf{A}$ from $x$ to $t$ such that the expected score decreases along that path. That is, we show that for all $n \in \{1, \ldots, N\}$, we have $\bar{S}(z^{(n-1)}, F) - \bar{S}(z^{(n)}, F) \geq 0$ and $\bar{S}(z^{(N)}, F) - \bar{S}(t, F) \geq 0$ with one inequality being strict under the conditions for strict $\mathcal{F}$-consistency. Then, a simple telescope sum argument concludes the proof.

The interpretation of condition 1 is obvious, meaning that the "start" of the path is at $x$. Condition 3 means that, for $N \geq 1$, the sequence approaches $t$ in the componentwise order. If the score satisfied some sort of componentwise order-sensitivity as discussed in Lambert et al. (2008) and in Section 3.2 of Fissler and Ziegel (2019) we would basically be done. Indeed, this componentwise order-sensitivity holds in the "quantile-components" 1 up to $k-1$. However, invoking Fissler and Ziegel (2019, Lemma B.2), componentwise order-sensitivity in all components would imply that $T$ would have elicitable components only, which is false. Even though componentwise order-sensitivity fails in the $k$th argument, one can show something related: For any $z, z' \in \mathsf{A}$ with $z_m = z'_m$ for all $m \in \{1, \ldots, k-1\}$ and $z_k \leq z'_k \leq -C(z, F)$ or $z_k \geq z'_k \geq -C(z, F)$ it holds that $\bar S(z, F) \geq \bar S(z', F)$ where

$$(6) \qquad C(z, F) := \sum_{m=1}^{k-1} \frac{p_m}{q_m} \left( \int_{(-\infty, z_m]} (z_m - y) \mathrm{d}F(y) - q_m z_m \right), \quad z \in \mathsf{A}.$$

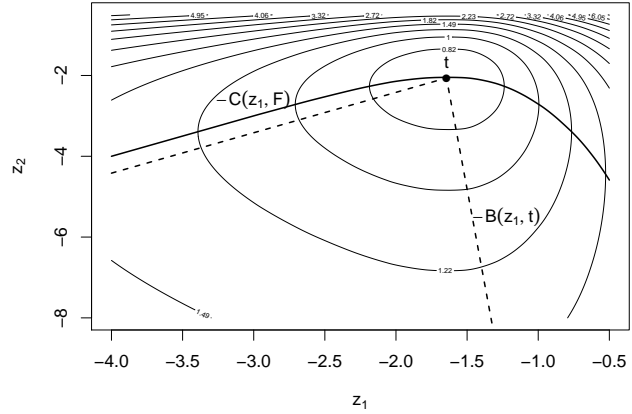

Figure 1: Graphical illustration of the function $z_1 \mapsto -C(z, F)$ (solid line) and $z_1 \mapsto -B(z, t)$ (dashed line) for $k = 2$, $F = \Phi$, a standard normal, and $t = \big(\mathrm{VaR}_{0.05}(F), \mathrm{ES}_{0.05}(F)\big) \approx (-1.64, -2.06)$. The contour lines correspond to the expected score $\bar S(z, F)$ for $G_1 = 0$, $\mathcal{G}_2(x) = -\log(-x)$, $x < 0$.

Figure 1 depicts the situation for $k = 2$ and $q_1 = 0.05$. Note that $C(z, F)$, $B(z, t)$ do not depend on $z_k$, so for $k = 2$ they can be displayed as a function of $z_1$. For illustration, we also display the level sets of the expected score $\bar S(z, F)$ with $G_1 = 0$, $\mathcal{G}_2(x) = -\log(-x)$, $x < 0$; compare Nolde and Ziegel (2017); Patton et al. (2019). As Figure 1 suggests, one can show that $t_k \geq -C(z, F)$. That means order-sensitivity in the $k$th argument may fail to hold when $x_k < t_k$. This is the reason why we impose the additional condition 4. It enforces that the expected score is first minimized in the quantile components until we are in a region where order sensitivity in the $k$th argument holds. To make this condition dependent on the particular distribution $F$ only via $t = T(F)$, we give the slightly stronger condition in terms of the function $B$ defined at (5). This exploits the fact that for any $F$ it holds that $-B(z, T(F)) \leq -C(z, F)$; see Figure 1. In summary, conditions 3 and 4 ensure that $\bar{S}(z^{(n-1)}, F) - \bar{S}(z^{(n)}, F) \geq 0$ for all $n \in \{1, \ldots, N\}$ if $N \geq 1$.

Finally, showing that $\bar{S}(z^{(N)}, F) - \bar{S}(t, F) \geq 0$ basically corresponds to the original proof given in Fissler and Ziegel (2016a). To show this final step, one needs to have condition 2, which was also not given in Fissler and Ziegel (2016a, Theorem 5.2). That means whenever $x \in \mathsf{A}$ and $t \in T(\mathcal{F})$ are such that $t_r, x_r \in \mathsf{A}'_{r,\min\{t_k,x_k\}}$ for all $r \in \{1, \ldots, k-1\}$, one can set $N = 0$. If this condition is not satisfied, then it is necessary to choose $N \geq 1$.

PROOF. Let $F \in \mathcal{F}$, $t = T(F)$ and $x \in \mathsf{A}$. Let $(z^{(n)})_{n=0,\ldots,N} \subseteq \mathsf{A}$ be a sequence that satisfies the above conditions. In the decomposition $\bar{S}(x, F) - \bar{S}(t, F) = R_1 + R_2$ on top of page 1702, we were not completely precise if $F$ is not continuous. For $a, b \in \mathbb{R}$, we define

$$I(a,b) = \begin{cases} (a,b], & \text{if } a \leq b \\ (b,a], & \text{if } a \geq b, \end{cases} \qquad \bar{I}(a,b) = \begin{cases} [a,b], & \text{if } a \leq b \\ [b,a], & \text{if } a \geq b. \end{cases}$$

For any $z, z' \in \mathsf{A}$ and any $w \in \mathsf{A}'_k$ it holds that $\bar{S}(z', F) - \bar{S}(z, F) = R_1 + R_2$ with

$$R_1 = \sum_{r=1}^{k-1} (F(z'_r) - q_r)\left(G_r(z'_r) + \frac{p_r}{q_r}G_k(w)z'_r\right)$$
$$- (F(z_r) - q_r)\left(G_r(z_r) + \frac{p_r}{q_r}G_k(w)z_r\right)$$
$$- \operatorname{sgn}(z'_r - z_r)\int_{I(z'_r, z_r)}\left(G_r(y) + \frac{p_r}{q_r}G_k(w)y\right)\mathrm{d}F(y).$$
$$R_2 = -\mathcal{G}_k(z'_k) + \mathcal{G}_k(z_k) + G_k(w)(z'_k - z_k) + (G_k(z'_k) - G_k(w))(z'_k + C(z', F))$$
$$- (G_k(z_k) - G_k(w))(z_k + C(z, F))$$

where $C(z, F)$ is defined at (6). If $t_r \leq z_r \leq z'_r$, and $y \mapsto G_r(y) + \frac{p_r}{q_r}G_k(w)y$ is increasing on $[z_r, z'_r]$ then the $r$th summand of $R_1$ is bounded below by

$$(F(z'_r) - q_r)\left(G_r(z'_r) + \frac{p_r}{q_r}G_k(w)z'_r\right) - (F(z_r) - q_r)\left(G_r(z_r) + \frac{p_r}{q_r}G_k(w)z_r\right)$$
$$- (F(z'_r) - F(z_r))\left(G_r(z'_r) + \frac{p_r}{q_r}G_k(w)z'_r\right)$$
$$= (F(z_r) - q_r)\left(G_r(z'_r) + \frac{p_r}{q_r}G_k(w)z'_r - G_r(z_r) - \frac{p_r}{q_r}G_k(w)z_r\right) \geq 0.$$

If $z'_r \leq z_r \leq t_r$, and $y \mapsto G_r(y) + \frac{p_r}{q_r}G_k(w)y$ is increasing on $[z'_r, z_r]$ then the $r$th summand of $R_1$ can be rewritten as

$$(F(z'_r) - q_r)\left(G_r(z'_r) + \frac{p_r}{q_r}G_k(w)z'_r\right) - (F(z_r) - q_r)\left(G_r(z_r) + \frac{p_r}{q_r}G_k(w)z_r\right)$$
$$+ \int_{[z'_r, z_r)}\left(G_r(y) + \frac{p_r}{q_r}G_k(w)y\right)dF(y)$$
$$- (F(z'_r) - F(z'_r-))\left(G_r(z'_r) + \frac{p_r}{q_r}G_k(w)z'_r\right)$$
$$+ (F(z_r) - F(z_r-))\left(G_r(z_r) + \frac{p_r}{q_r}G_k(w)z_r\right)$$
$$\geq (F(z'_r-) - q_r)\left(G_r(z'_r) + \frac{p_r}{q_r}G_k(w)z'_r\right)$$
$$- (F(z_r-) - q_r)\left(G_r(z_r) + \frac{p_r}{q_r}G_k(w)z_r\right)$$
$$+ (F(z_r-) - F(z'_r-))\left(G_r(z'_r) + \frac{p_r}{q_r}G_k(w)z'_r\right)$$
$$= (F(z_r-) - q_r)\left(G_r(z'_r) + \frac{p_r}{q_r}G_k(w)z'_r - G_r(z_r) - \frac{p_r}{q_r}G_k(w)z_r\right) \geq 0.$$

In case that $z'_r < z_r = t_r$ or $z'_r > z_r = t_r$, if $y \mapsto G_r(y) + \frac{p_r}{q_r}G_k(w)y$ is strictly increasing on $\bar{I}(z'_r, t_r)$, and $F$ has a unique $q_r$-quantile, then even $R_1 > 0$. Indeed, if $F$ does not put any mass on the open interval $\text{int}(I(z'_r, t_r))$, then the last inequality is strict in both cases. Otherwise, the first inequality is strict because $y \mapsto G_r(y) + \frac{p_r}{q_r}G_k(w)y$ is strictly increasing. In summary, setting $w = \min\{z'_k, z_k\}$, exploiting condition 2 and the convexity of $\mathsf{A}$, this implies that both $R_1 \geq 0$ for $z' = z^{(n-1)}$, $z = z^{(n)}$, $n \in \{1, \ldots, N\}$, and also for $z' = z^{(N)}$, $z = t$. For the latter case, we even have $R_1 > 0$ under the conditions for strict $\mathcal{F}$-consistency and if $(z_m^{(N)})_{m=1,\ldots,k-1} \neq (t_m)_{m=1,\ldots,k-1}$.

Discussing $R_2$, we use that $-t_k = \sum_{m=1}^{k-1}(p_m/q_m)(\int_{(-\infty,t_m]}(t_m-y)dF(y) - q_m t_m)$ to check that for any $z \in \mathsf{A}$

$$(7) \qquad -t_k \leq C(z, F) \leq B(z, t),$$

where both inequalities are equalities for $z = t$. For the first inequality, note that the summands of $-t_k$ are the expectations of the consistent $q_m$-quantile scoring functions $\mathbb{1}\{y \leq w\}(w - y) - q_m w$ (up to scaling with $p_m/q_m$) evaluated at the quantile $t_m$. On the other hand, the summands of $C(z, F)$ are the expectations of the same scoring functions evaluated at arbitrary other points $z_m$, and therefore the first inequality holds due to the consistency of the quantile score. For the second inequality, we replace $-t_k$ in the definition of $B(z,t)$ with the given formula. This yields

$$B(z,t) = \sum_{m=1}^{k-1} \frac{p_m}{q_m} \Big( \int_{(-\infty,t_m]} (t_m - y) \mathrm{d}F(y) - z_m q_m - (t_m - z_m) \mathbb{1}\{t_m < z_m\} \Big).$$

The inequality follows by considering for each summand the cases $t_m < z_m$ and $t_m \geq z_m$ separately, using that

$$\int_{(-\infty,t_m]} (t_m - y)\mathrm{d}F(y) = \int_{(-\infty,z_m]} (z_m - y)\mathrm{d}F(y) - F(z_m)(z_m - t_m)$$
$$+ \mathrm{sgn}(t_m - z_m) \int_{I(t_m,z_m)} (t_m - y)\mathrm{d}F(y).$$

Let $z' = z^{(N)}$, $z = t$ and again $w = \min\{z'_k, z_k\}$. Then due to (7) and the monotonicity of $G_k$ the term $R_2$ is bounded below by

$$-\mathcal{G}_k(z'_k) + \mathcal{G}_k(t_k) + G_k(w)(z'_k - t_k) + (G_k(z'_k) - G_k(w))(z'_k - t_k) \geq 0$$

invoking the convexity of $\mathcal{G}_k$. The inequality becomes strict if $\mathcal{G}_k$ is strictly convex and $z'_k \neq t_k$. Now, let $z' = z^{(n-1)}$ and $z = z^{(n)}$ for $n \in \{1,\ldots,N\}$ and again $w = \min\{z'_k, z_k\}$. If $t_k \leq z_k \leq z'_k$, then we obtain

$$R_2 \geq -\mathcal{G}_k(z'_k) + \mathcal{G}_k(z_k) + G_k(z'_k)z'_k - G_k(z_k)z_k + (G_k(z'_k) - G_k(z_k))(-t_k)$$
$$= -\mathcal{G}_k(z'_k) + \mathcal{G}_k(z_k) + G_k(z'_k)(z'_k - t_k) - G_k(z_k)(z_k - t_k)$$
$$\geq -\mathcal{G}_k(z'_k) + \mathcal{G}_k(z_k) + G_k(z'_k)(z'_k - t_k) - G_k(z'_k)(z_k - t_k) \geq 0,$$

where the penultimate inequality is due to the fact that $G_k$ is increasing and the last inequality follows due to the convexity of $\mathcal{G}_k$. The last inequality is strict if $\mathcal{G}_k$ is strictly convex and if $z_k \neq z'_k$. If $z'_k \leq z_k \leq t_k$ and $z_k \leq -B(z,t) = -B(z',t)$, then, using (7)

$$R_2 = -\mathcal{G}_k(z'_k) + \mathcal{G}_k(z_k) + G_k(z'_k)z'_k - G_k(z_k)z_k + (G_k(z'_k) - G_k(z_k))C(z,F)$$
$$\geq -\mathcal{G}_k(z'_k) + \mathcal{G}_k(z_k) + G_k(z'_k)z'_k - G_k(z_k)z_k + (G_k(z'_k) - G_k(z_k))(-z_k)$$
$$\geq 0$$

by convexity of $\mathcal{G}_k$. Again, the inequality is strict if $\mathcal{G}_k$ is strictly convex and if $z_k \neq z'_k$. In summary, $R_2 \geq 0$ for $z' = z^{(n-1)}$, $z = z^{(n)}$, $n \in \{1,\ldots,N\}$, and also for $z' = z^{(N)}$, $z = t$. For the latter case, we even have $R_2 > 0$ if $\mathcal{G}_k$ is strictly convex and if $t_k \neq z_k^{(N)}$. $\square$

We would like to remark that the additional condition stated in Proposition 2 holds in most practically relevant cases of action domains $\mathsf{A}$. In fact, in the following situations, one can set $N = 0$ meaning that the original proof in Fissler and Ziegel (2016a) is applicable:

- $\mathsf{A} = \mathbb{R}^k$.
- $\mathsf{A} = \mathsf{A}_0$, where $\mathsf{A}_0$ is the maximal sensible action domain for the functional $T$ defined in Theorem 5.2 in Fissler and Ziegel (2016a) for $q_1 < \cdots < q_{k-1}$.[1] It is given by

$$\mathsf{A}_0 := \Big\{x \in \mathbb{R}^k : x_1 \leq \cdots \leq x_{k-1}, \ x_k \leq \sum_{m=1}^{k-1} p_m x_m \Big\}.$$

  In particular, this construction retrieves the result of Corollary 5.5 in Fissler and Ziegel (2016a), considering the maximal action domain $\mathsf{A}_0 = \{x \in \mathbb{R}^2 : x_1 \geq x_2\}$ for the functional $T = (\mathrm{VaR}_\alpha, \mathrm{ES}_\alpha)$.
- Also, $\mathsf{A} = \mathsf{A}_0^+ := \mathsf{A}_0 \cap \mathbb{R}_+^k$ and $\mathsf{A} = \mathsf{A}_0^- := \mathsf{A}_0 \cap \mathbb{R}_-^k$ satisfy the condition.
- The action domain considered in Nolde and Ziegel (2017, Theorem C.3), corresponding to $\mathsf{A} = \mathbb{R} \times (-\infty, 0)$ in our sign convention, also satisfies the condition given in Proposition 2.

For the action domain $\mathsf{A} = \{x \in \mathbb{R}^2 : x_2 \leq x_1 \leq x_2 + c\}$ considered in Fissler and Ziegel (2019, Proposition 4.10) one generally needs that $N > 0$. However, the existence and construction of the sequence $(z^{(n)})_{n=0,\ldots,N}$ is obvious.

Acerbi and Szekely (2014) introduced a family of scoring functions $S^W$, $W \in \mathbb{R}$, with corresponding action domains $\mathsf{A}^W = \{x \in \mathbb{R}^2 : x_2 > Wx_1\}$; see page 1697 in Fissler and Ziegel (2016a) for a discussion. The results are as follows:

- For $W > 1$ the additional condition of Proposition 2 is satisfied (note that the potentially problematic point $(0,0)$ is not in $\mathsf{A}^W$ then.)
- For $W = 1$, the condition is empty since $\mathsf{A}^1 \cap T(\mathcal{F}) \subseteq \mathsf{A}^1 \cap \mathsf{A}_0 = \emptyset$.
- For $W \in (0, 1)$ the condition fails to be satisfied.
- For $W = 0$ the condition holds.
- For $W \in [(\alpha - 1)/\alpha, 0)$ the condition fails to be satisfied.
- For $W < (\alpha - 1)/\alpha$ the condition holds.

**3. Remarks on Propositions 4.2 and 4.4.** In Fissler and Ziegel (2016a) Proposition 4.2(i) we claimed that $g_m > 0$ and in Proposition 4.4(i) we claimed that the matrix-valued function $h = (h_{r,l})_{r,l=1,\ldots,k}$ is positive definite on $\mathrm{int}(\mathsf{A})$. In the proof of Theorem 5.2(iii), we made similar claims about the function $g_k$ and the derivatives of the functions appearing at (5.3). It turns out that we were slightly imprecise with these claims. Indeed, these functions are all non-negative (positive semi-definite). Moreover, we can show that the functions $g_m$ in Proposition 4.2(i) (as well as the functions appearing in the proof of Theorem 5.2(iii)) are strictly positive almost everywhere.

---

[1] With 'maximal sensible' we mean that $T(\mathcal{F}) \subseteq \mathsf{A}_0$ for any class of distributions $\mathcal{F}$.

We will first demonstrate that the function $h$ appearing in Proposition 4.4(i) is positive semi-definite. Assume there is some $t \in \mathrm{int}(\mathsf{A})$ and some $v \in \mathbb{S}^{k-1}$ such that $v^\top h(t) v < 0$. Recall that due to the assumptions and Theorem 3.2, the function $h$ is continuous. Hence, there is an open neighbourhood $t \in U \subseteq \mathrm{int}(\mathsf{A})$ such that $v^\top h(x) v < 0$ for all $x \in U$. Invoking the surjectivity of $T$, let $F \in \mathcal{F}$ such that $T(F) = t$. Then there is an $\varepsilon > 0$ such that – using the notation on top of page 1701 in Fissler and Ziegel (2016a)

$$\frac{\mathrm{d}}{\mathrm{d}s}\bar{S}(t+sv, F) = \bar{q}(F)sv^\top h(t+sv)v \begin{cases} > 0, & \text{for } s \in (-\varepsilon, 0), \\ < 0, & \text{for } s \in (0, \varepsilon). \end{cases}$$

This is a contradiction to the $\mathcal{F}$-consistency of $S$. Therefore, the function $\phi$ with Hessian $h$ is convex (but not necessarily strictly convex).

Now, we will show that the functions $g_m$ appearing in Proposition 4.2(i) are positive almost everywhere, the arguments for the remaining functions being similar. Let $m \in \{1, \ldots, k\}$. The argument at the top of page 1701 in Fissler and Ziegel (2016a) show that for all $x \in \mathrm{int}(\mathsf{A}'_m)$ there is an $\varepsilon > 0$ such that $g_m > 0$ on $(x - \varepsilon, x + \varepsilon) \setminus \{x\}$. Due to a compactness argument, $g_m$ is positive on any compact set almost everywhere. Since $\mathrm{int}(\mathsf{A}'_m)$ is $\sigma$-compact, $g_m > 0$ almost everywhere. The continuity of $g_m$ also implies that $g_m \geq 0$ everywhere.

## References.

Institute for Statistics and Mathematics
WU Vienna University of Economics and Business
Welthandelsplatz 1, 1020 Vienna
Austria
E-mail: tobias.fissler@wu.ac.at

Institute of Mathematical Statistics
 and Actuarial Science
University of Bern
Alpeneggstrasse 22, 3012 Bern
Switzerland
E-mail: johanna.ziegel@stat.unibe.ch